\documentclass[final,onefignum,onetabnum]{siamart190516}

\usepackage{lipsum}
\usepackage{amsfonts}
\usepackage{graphicx}
\usepackage{epstopdf}
\usepackage{tikz}
\usepackage[caption=false]{subfig}
\usepackage{algorithmic}
\ifpdf
  \DeclareGraphicsExtensions{.eps,.pdf,.png,.jpg}
\else
  \DeclareGraphicsExtensions{.eps}
\fi

\headers{Infinite domain Helmholtz equation}{Anil Zengino\u{g}lu}

\usepackage{mathrsfs}

\DeclareSymbolFontAlphabet{\mathrsfs}{rsfs}

\newcommand{\be}{\begin{equation}}
\newcommand{\ee}{\end{equation}}
\newcommand{\iu}{{i\mkern1mu}}

\title{A null infinity layer for wave scattering}

\author{Anil Zengino\u{g}lu\thanks{Institute for Physical Science and Technology, University of Maryland, College Park, MD 20742
  (\email{anil@umd.edu)}.}}

\usepackage{amsopn}

\date{}

\begin{document}

\maketitle

\begin{abstract}
We show how to solve time-harmonic wave scattering problems on unbounded domains without truncation. The technique, first developed in numerical relativity for time-domain wave equations, maps the unbounded domain to a bounded domain and scales out the known oscillatory decay towards infinity. We design a null infinity layer that corresponds to the infinite exterior domain and restricts the transformations to an annular domain. The method does not require the local Green function. Therefore we can use it to solve Helmholtz equations with variable coefficients and certain nonlinear source terms. The method's main advantages are the exact treatment of the local boundary and access to radiative fields at infinity. The freedom in the transformations allows us to choose parameters adapted to high-frequency wave propagation in the exterior domain. We demonstrate the efficiency of the technique in one- and two-dimensional numerical examples.
\end{abstract}



\section{Introduction}

We consider time-harmonic wave scattering from a bounded obstacle $D \subset \mathbb{R}^d$ in $d$ dimensions. The scattered wave, $U$, satisfies the Helmholtz equation and the Sommerfeld radiation condition at infinity
\begin{equation}
	\label{eq:helm}
	\begin{gathered}
		\Delta U + k^2 U = F \quad \mathrm{in } \ \mathbb{R}^d\setminus \bar{D}, \\
		\partial_r U - \iu k U = o\left(r^{\frac{1-d}{2}}\right) \quad \mathrm{as} \ r\to\infty,\ d=1,2,3.
	\end{gathered}
\end{equation}
Dirichlet or Neumann conditions are given on the surface of the scatterer $D$.

Two essential difficulties for this problem are: (i) the unbounded domain and (ii) the highly oscillatory behavior of the solution for high wave numbers $k$. Both difficulties have been active areas of research with an extensive list of proposed treatments (see reviews \cite{TSYNKOV1998465, HagstromReview, engquist2003computational}).

Methods for handling the unbounded domain fall into two main categories: local and exact \cite{TSYNKOV1998465, HagstromReview, givoli2013numerical}. Local methods are convenient but approximate, while exact methods are accurate but cumbersome. Local methods are most common due to their convenience and generality. They include high-order boundary conditions \cite{EngquistMajda77, BaylissTurkel80} and absorbing regions \cite{israeli1981approximation, BERENGER1994185}. There are spectrally convergent local methods that may be sufficiently accurate for many applications \cite{hagstrom2009complete}. Exact methods, on the other hand, provide better accuracy but are non-local and therefore computationally expensive \cite{keller1989exact, GROTE1998327, grote1996nonreflecting}. Sophisticated methods exist for the fast evaluation of the associated kernels \cite{alpert2000rapid}, but they are hard to implement. In summary, there are excellent methods available to deal with outer boundaries that achieve great accuracy for certain sets of problems. There is, however, no {\it local and exact} treatment of the numerical outer boundary, and the topic is an active area of research \cite{9196168, kirby2020finite, papadimitropoulos2020double, petropavlovsky2020numerical, duhamel2020computation}.

In principle, a local and exact method to compute unbounded domain solutions to Helmholtz equations would be compactification. Mapping the unbounded solution domain to a bounded numerical domain would allow us to solve the equation without the need for an outer boundary treatment \cite{GroschOrszag77, boyd1982optimization, shen2009some}. However, this method leads to slow convergence for wave equations due to infinite oscillations in the asymptotic region \cite{GroschOrszag77, shen2014approximations}. The mapped finite grid cannot represent infinite oscillations.

In the time-domain, a suitable transformation of the time coordinate solves the infinite resolution problem for the wave equation \cite{ZENGINOGLU20112286}. Even though there is no time coordinate in the Helmholtz equation, one can still perform the time transformation as a rescaling of the unknown variable that takes out the known oscillatory behavior of the solution in the asymptotic domain \cite{ZengFramework}.

The transformation is inspired by Penrose compactification and the notion of conformal infinity in relativity \cite{Penrose, Penrose65, Frauendiener2004, Winicour2012}. It consists of the following three steps:
\begin{itemize}
	\item[(i)] Scale-out the asymptotic fall-off behavior.
	\item[(ii)] Scale-out the known oscillatory behavior.
	\item[(iii)] Map the unbounded domain to a bounded domain.
\end{itemize}
This method is an application of hyperboloidal compactification (\cite{Zenginoglu08, ZENGINOGLU20112286, Macedo_2020}) to time-harmonic solutions \cite{ZengFramework, ansorg2016spectral, macedo2018hyperboloidal, bizon2020toy, jaramillo2021pseudospectrum, jaramillo2021gravitational, destounis2021pseudospectrum, gasperin2021physical}. It includes compactification along characteristic surfaces as a special case. The outer boundary of the mapped solution domain corresponds to null infinity. Therefore, we refer to the method as \emph{null infinity compactification} (NIC).

This paper uses null infinity compactification to solve the Helmholtz equation on unbounded domains with local methods. Access to the entire solution leads to an exact method that does not require an approximate boundary treatment and has only discretization errors. The lack of a numerical boundary treatment improves computational efficiency. In addition, one has direct access to the solution at infinity, simplifying the extraction of radiative solutions in the far-field. Such access is essential in many applications where limits at infinity provide observables of interest, such as the echo area in acoustics, the radar cross-section in electromagnetics, or the gravitational waveform in general relativity.

Steps (i) and (ii) transform the generic behavior of solutions to the Helmholtz equation from algebraic decay with oscillations towards infinity to no decay without oscillations towards infinity. In geometric terms, step (i) corresponds to a conformal rescaling, and step (ii) corresponds to a time transformation in frequency domain \cite{ZengFramework}. This separation of the unknown into oscillatory and non-oscillatory parts is similar to the geometrical optics decomposition into amplitude, and phase \cite{engquist2003computational}. The difference here is that the phase is now a prescribed function that takes care of asymptotic oscillations. 

A significant advantage of null infinity compactification in the frequency domain compared to the time domain is the adaptation to high-frequency wave propagation. In the time domain, the Courant-Friedrichs-Lewy condition restricts the transformation \cite{ZENGINOGLU20112286}. In the frequency domain, no such restriction exists, and we can adapt the transformation to high-frequency wave propagation problems. 

Restricting the transformations to an annular domain gives us a layer corresponding to the infinite exterior domain with the boundary at null infinity. This null infinity layer (NIL) differs from the hyperboloidal layer of \cite{ZENGINOGLU20112286} because it can use both characteristic and hyperboloidal coordinates. NIL with characteristic coordinates is similar to the perfectly absorbing layer of \cite{wang2017perfect, yang2021truly} except that the solution is not artificially damped in the layer, and the radiative solution at infinity can be obtained at the numerical outer boundary. In addition, NIL does not rely on the local Green function. The method can be applied to equations with variable coefficients or nonlinear source terms if certain asymptotic fall-off conditions are satisfied. 

In the next section, we present NIC on a simple one-dimensional example. In Sec.~\ref{sec:nic}, we specify the proposed transformation for the constant-coefficient Helmholtz equation. In Sec.~\ref{sec:nil} we present the restriction of the transformations to a layer. Numerical experiments are presented in Sec.~\ref{sec:numerical}. 

\section{Motivation for the transformation: A simple example}\label{sec:simple}
Consider the Helmholtz equation in one dimension with the Sommerfeld boundary condition 
\begin{align}  
\label {eq:Helmholtz1d} d_r^2 U + k^2 U &= 0, \quad r\in[0,\infty),\\ 
\label{eq:sommerfeld1d}d_r U - \iu k U &= 0 \quad \mathrm{as} \ r\to\infty.
\end{align}
The Helmholtz equation has oscillatory solutions. For example, the plane wave $U(r)=e^{\iu k r}$ solves the above system.

Mapping the infinite domain to a finite domain leads to an equation with a singular term. We map the unbounded domain, $r\in[0,\infty)$, to a compact domain adding the point at infinity, $\rho\in[0,\tfrac{\pi}{2}]$, using the transformation
\be\label{eq:space_compact1d} r = \tan \rho \ee
We get for \eqref{eq:Helmholtz1d}
\be\label{eq:simple} d_\rho(\cos^2 \rho \ d_\rho U )+ \frac{k^2}{\cos^2 \rho} U = 0. \ee
The transformed wavenumber is unbounded near the boundary at $\rho=\tfrac{\pi}{2}$. When we map an infinite domain to a finite domain, the constant wavenumber of oscillatory solutions maps to an unbounded wavenumber near the new domain boundary. This behavior is the reason that compactification is considered ineffective for equations with oscillatory solutions \cite{GroschOrszag77, shen2014approximations}. We cannot resolve the plane wave solution in the new coordinate, $U(\rho)=e^{\iu k \tan\rho}$, near the domain boundary.  

However, we know the asymptotic form of these oscillations so we can scale them out. Define
\be\label{eq:trafo1d} U(\rho) = e^{\iu k (\tan\rho - \rho)} u(\rho). \ee
The rescaled equation becomes
\be\label{eq:transformed} d_\rho(\cos^2 \rho \ d_\rho u ) + 2\iu k \sin^2\rho \ d_\rho u + \left(k^2 (1+\sin^2\rho)+\iu k \sin(2\rho) \right) u = 0. \ee
This equation is well-behaved at the boundary as opposed to \eqref{eq:simple}. The rescaling transforms the plane wave in $r$ to a plane wave in $\rho$ that reads $u(\rho)=e^{\iu k \rho}$.

Evaluating the equation at $\rho=\tfrac{\pi}{2}$, we get the Sommerfeld condition for the transformed variable
\[ d_\rho u - \iu k u = 0. \]
The condition is not supplied separately. It is a consequence of the equation evaluated at the outer boundary. The boundary condition is behavioral, so no boundary treatment is necessary. 

\subsection{Relation to time transformations}\label{sec:time}
We recap the derivation of the Helmholtz equation from the wave equation to connect the rescaling \eqref{eq:trafo1d} to time transformations. Consider the scalar wave equation with unit speed of propagation 
\[ -\partial_t^2 V + \partial_r^2 V = 0\,, \]
We look for solutions with a single frequency, $k$. The ansatz 
\be\label{eq:ansatz} 
V(r,t)= e^{-\iu k t} U(r),
\ee 
leads to the Helmholtz equation \eqref{eq:Helmholtz1d}.

We use the spatial mapping \eqref{eq:space_compact1d} for null infinity compactification in time-domain. To avoid infinite oscillations at the domain boundary, a new time coordinate, $\tau$, must be introduced. There is large freedom in the choice of coordinates. We require that the new time coordinate satisfies $\partial_\tau=\partial_t$ so that the transformed equation has time-independent coefficients. In relativity, this requirement corresponds to the invariance of the time translation symmetry of the underlying Minkowski metric. The time transformation takes the form
\be\label{eq:tau} \tau = t - h(\rho),\ee
where $h(\rho)$ is the height function. We make the single frequency ansatz
\[ V(\rho, \tau) = e^{-\iu k \tau} u(\rho) = e^{-\iu k t} e^{\iu k h(\rho)} u(\rho)\,.\]
Comparing with \eqref{eq:ansatz}, we see that the time transformation \eqref{eq:tau} in time domain corresponds to a rescaling, or more accurately, a phase shift $U(\rho) = e^{\iu k h(\rho)} u(\rho)$ in frequency domain \cite{ZengFramework, marchner2021stable}.

The height function in this example, $h(\rho)=\tan\rho - \rho$, is such that the coordinates satisfy the following convenient relationship
\be\label{eq:time_chars} \tau - \rho = t - r\,.\ee
Outgoing characteristics in standard coordinates $(t,r)$ have the same form as outgoing characteristics in the transformed coordinates $(\tau,\rho)$. One can derive useful prescriptions for the height function by imposing conditions on the form of outgoing characteristics in compactifying coordinates \cite{ZENGINOGLU20112286, bernuzzi2011binary}.

\subsection{Dispersion relation}\label{sec:dispersion}
We derive dispersion relations for the three variants of the Helmholtz equation using a plane wave ansatz. With $U = e^{\iu \xi r}$ for \eqref{eq:Helmholtz1d}, we get the usual dispersion relation $\xi_\pm=\pm k$ corresponding to outgoing and incoming waves. The Sommerfeld condition \eqref{eq:sommerfeld1d} breaks the symmetry between outgoing and incoming waves by enforcing that there are no incoming waves from infinity. Only outgoing waves with $\xi = k$ survive asymptotically. 

The same symmetry between outgoing and incoming waves is present in the compactified version \eqref{eq:simple}. Using the ansatz $U=e^{\iu \xi \rho}$, the dispersion relation reads
\[ \cos^2\rho\ \xi^2 + 2\iu \cos\rho \sin\rho\ \xi - \frac{k^2}{\cos^2\rho} = 0, \]
with solutions
\be\label{eq:divergent_dispersion} \xi_\pm = \pm \sqrt{\frac{k^2}{\cos^4\rho}-\tan^2\rho} \mp \iu \tan\rho. \ee
The dispersion relation in the spatially compactified coordinates shows that both outgoing and incoming wave numbers blow up and cannot be resolved near the domain boundary at $\rho=\pi/2$.  

The dispersion relation for the transformed Helmholtz equation \eqref{eq:transformed} with the ansatz $u=e^{\iu k \rho}$ reads
\[ (k-\xi)^2 \sin^2\rho + \iu \sin(2\rho) (k-\xi) + (k^2-\xi^2) = 0. \]
A solution to this relation is $\xi_+=k$, just as the standard Helmholtz equation for outgoing waves. This relationship arises by construction of the time transformation in \eqref{eq:time_chars}.
The solution that corresponds to incoming waves is singular at the outer boundary
\be\label{eq:incoming} \xi_-=-\frac{k(1+\sin^2\rho)+2\iu \sin(2\rho)}{\cos^2\rho}. \ee
The divergence near the boundary is stronger than for \eqref{eq:divergent_dispersion}. The incoming wavenumber is larger than the outgoing wavenumber throughout the domain, $|\xi_-|>|\xi_+|$ for all $\rho>0$. Adapting the coordinates to outgoing waves decreases the resolution of incoming waves in the domain \cite{calabrese2006asymptotically}. 

Lower resolution for incoming waves is undesirable near scatterers of arbitrary shape because waves propagate in all directions. In higher dimensions, the restriction of the transformations to a layer may be helpful. Next, we present the general transformations, and in Sec.~\ref{sec:nil} we discuss the null infinity layer, which restricts these transformations to a layer.

\section{Null infinity compactification for the Helmholtz equation}\label{sec:nic}
The simple example presented in the previous section generalizes to arbitrary dimensions with a few modifications. In higher dimensions, we remove the decay of the oscillatory solution to ensure regularity of the transformed equation. The transformation is performed along the outgoing direction towards infinity. We summarize the transformation with the following directive: Scale-out the oscillatory decay and compactify
\be\label{eq:summary} u(\rho, \omega^{d-1}) = r(\rho)^{\tfrac{d-1}{2}} e^{-\iu k h(\rho)} U(r(\rho), \omega^{d-1})\,.\ee
Here, $\omega^{d-1}$ are angular coordinates on a $d-1$ dimensional sphere, $h(\rho)$ is a suitably chosen height function satisfying certain asymptotic conditions, and $r(\rho)$ is the compactification of the radial coordinate $r$. Below we present these steps in detail on the Helmholtz equation in spherical coordinates
\begin{equation}\label{eq:helm_sph}
	\partial_r^2 U + \frac{d-1}{r} \partial_r U + \frac{1}{r^2}\triangle_{S^{d-1}} U + k^2 U = F,
\end{equation}
where $\triangle_{S^{d-1}}$ is the Laplace operator on the $d-1$ dimensional sphere $S^{d-1}$. We perform the steps listed in the Introduction to derive the null infinity compactification of this equation.

(i) We introduce the rescaled variable $\tilde{U}$ which scales out the fall-off behavior of the unknown near infinity via $\tilde{U} = r^{\tfrac{d-1}{2}} U$. The rescaling removes the first order derivative in $r$ in exchange for a low-order term that falls off as $r^2$. The Helmholtz equation becomes after division by $r^{\tfrac{1-d}{2}}$
\[ \partial_r^2 \tilde{U} + \frac{1}{r^2}\triangle_{S^{d-1}} \tilde{U} + \left(k^2 - \frac{(1-d)(3-d)}{4 r^2} \right) \tilde{U} = F r^{\tfrac{d-1}{2}}. \]
We will see in step (iii) that the $r^{-2}$ behavior of the terms that vanish for $d=1$ is just right for the regularity of the transformed equation when $d>1$. 

With this rescaling, the value of the unknown at infinity is a non-vanishing constant. The Sommerfeld radiation condition becomes
\[ \partial_r \tilde{U} - \iu k \tilde{U} = o(1) \quad \mathrm{as} \ r\to\infty. \]

(ii) In the second step, we introduce a rescaled variable $u$ which scales out the oscillatory behavior of the unknown near infinity via $u=e^{-\iu k h(r)} \tilde{U}$. This rescaling is the crucial step in Eq.~\eqref{eq:trafo1d}, and is related to a transformation of the time coordinate in the wave equation \cite{ZengFramework} as we demonstrated on the simple example in Sec.~\ref{sec:time}. We refer to $h(r)$ as the height function as is common in relativity \cite{reinhart1973maximal, beig1996vacuum}. The Helmholtz equation becomes after division by $e^{\iu k h}$
\begin{align*}
	\partial_r^2 u & + 2 \iu k H \partial_r u + \frac{1}{r^2}\triangle_{S^{d-1}} u +                                        \\
	        & \left[k^2 \left(1- H^2\right) - \frac{(1-d)(3-d)}{4 r^2} + \iu k \, d_r H \right] u = F r^{\tfrac{d-1}{2} e^{-\iu k h}},
\end{align*}
where $H:=d_r h$ is the radial derivative of the height function and is called the boost in analogy to Lorentz boosts. Inspecting the terms in the equation with respect to the powers of $k$, we recognize the analogs of eikonal and transport equations. The height function, $h$, plays the role of the phase, and the unknown, $u$, plays the role of the amplitude. The difference is that here the height function is prescribed explicitly such that the boost, $H$, satisfies \cite{Zenginoglu08, ZengFramework, jaramillo2021pseudospectrum}
\begin{equation} \label{eq:h_cond}
	H \leq 1, \quad \lim_{r\to\infty} H = 1, \quad \lim_{r\to\infty} d_r H = 0. 
\end{equation}
The Sommerfeld radiation condition written in terms of the rescaled unknown $u$ is independent of the wave number
\[ \partial_r u = o(1) \quad \mathrm{as} \ r\to\infty. \]

The rescaling removes asymptotic oscillations from the solution, which makes the equation amenable to compactification. There is considerable freedom in specific choices for the height function that satisfy \eqref{eq:h_cond} (see \cite{Zenginoglu08} for asymptotic conditions in asymptotically flat spacetimes and \cite{Macedo_2020} for a review of choices in the context of rotating black hole spacetimes). 

For example, setting $h(r)=r$, with $H=1$, removes all oscillations for an outgoing spherical wave centered at the origin. This choice corresponds to using the outgoing characteristics, $t-r$, as a time coordinate. The outgoing spherical wave in three dimensions, $U(r)=e^{\iu k r}/r$, becomes $u=1$. Various methods use similar rescalings to study Helmholtz equations. The infinite element method \cite{demkowicz2006few} and the perfectly absorbing layer \cite{wang2017perfect, yang2021truly} both use this rescaling to remove the oscillations from the outgoing solution. The pole condition in \cite{schmidt2008pole} is also related to this rescaling through the Laplace-transformation in combination with the scaling out of the asymptotic decay. The height function approach generalizes these methods in a geometric framework and increases the flexibility of null infinity compactification for handling heterogeneous media and smoothly matched layers. When the boost is strictly less than unity, we get compactification along hyperboloidal time surfaces \cite{Zenginoglu08}. 

(iii) Having removed the asymptotic decay and oscillations from the solution, we map the infinite domain in $r$ to a finite domain in a new radial coordinate $\rho$. Such mappings are well-known and have been extensively studied \cite{GroschOrszag77, boyd1982optimization, shen2009some, wang2017perfect}. Consider a mapping of the form
\[ r = g(\rho), \quad r\in[0,\infty),\ \rho\in[0,S), \]
such that
\begin{align*}
	 & \frac{dr}{d\rho} = g'(\rho) \equiv \frac{1}{G(\rho)} > 0, \quad \rho\in [0,S), \\
	 & g(0)      = 0, \qquad g(S) = \infty.
\end{align*}
Here, $S$ is the location of the outer boundary. The Helmholtz equation becomes after a division by $G(\rho)$
\begin{align}
	\label{eq:hyp_general}
	G \partial_\rho^2 u & + \left( d_\rho G + 2 \iu k H \right) \partial_\rho u + \frac{1}{G g^2}\triangle_{S^{d-1}} u +  \nonumber                        \\
	               & \left[k^2 \frac{1-H^2}{G} - \frac{(1-d)(3-d)}{4} \frac{1}{G g^2} + \iu k d_\rho H \right] u = \frac{F}{G} g^{\tfrac{d-1}{2}} e^{-\iu k h}.
\end{align}
This equation is the general form of the Helmholtz equation under the proposed null infinity compactification. We recover the standard Helmholtz equation for $g(\rho) = \rho$ implying $G(\rho)=1$ and $H(\rho)=0$. The regularity of the equation depends on the source function $F$. In our experiments, we consider source-free equations, but the method can incorporate compactly supported or sufficiently fast decaying source terms. The boost conditions and compactification guarantee the regularity of the other terms with divisions by $G$. We present a few specific examples next.

\subsection{Specific choices of free functions}
There is considerable freedom in the mapping $g(\rho)$ and the height function $h(\rho)$. A numerical analysis of the various choices is outstanding. The optimal choice will depend on the problem, but generally, the regularity of the transformed Helmholtz equation \eqref{eq:hyp_general} follows if $1-H^2\sim G\sim g^{-2}$ near the outer boundary at $\rho=S$. 

The simple example presented in Sect.~\ref{sec:simple} includes the following choices
\[ d=1, \quad F(\rho)=0, \quad g(\rho)=\tan \rho, \quad G(\rho) = \cos^2\rho, \quad h(\rho) =\tan\rho - \rho, \quad H(\rho)= \sin^2\rho\,.\]
The outer boundary is at $S=\tfrac{\pi}{2}$. We have $(1-H^2)/G = 1+\sin^2\rho$, and $G g^2 = \sin\rho$, so the transformed equation \eqref{eq:transformed} is regular at the outer boundary. 

For most numerical experiments presented in this paper, we specify the mapping $g$ and the height function $h$ as follows
\begin{align} 
\label{eq:g}
g(\rho) = \frac{\rho}{1-\rho} \quad &\mathrm{with}\quad G(\rho) \equiv \left(\frac{d g(\rho)}{d\rho}\right)^{-1} = (1-\rho)^2\,, \\
\label{eq:h}
h(\rho) = g(\rho) - \frac{\rho}{K} \quad &\mathrm{with}\quad H(\rho) \equiv G(\rho) \frac{d h(\rho)}{d\rho} = 1 - \frac{\Omega(\rho)}{K} 
\end{align}
where $K>0$ is a parameter that will be used to adapt the transformation to high-frequency wave propagation, and we define $\Omega(\rho):=\sqrt{G(\rho)} = 1-\rho$. 
The resulting Helmholtz equation reads
\begin{align}\label{eq:hyp_specific}
	\Omega^2 \partial_\rho^2 u & - 2\left(\Omega - \iu k \left(1-\frac{\Omega^2}{K}\right)\right) \partial_\rho u + \frac{1}{\rho^2}\triangle_{S^{d-1}} u + \nonumber \\
	             & \left[k^2 \left( \frac{2}{K} - \frac{\Omega^2}{K^2} \right) - \frac{(1-d)(3-d)}{4 \rho^2} + 2 \iu k \frac{\Omega}{K} \right] u = \frac{F}{\Omega^2} \left(\frac{\rho}{\Omega}\right)^{\tfrac{d-1}{2}} e^{-\iu k h}
\end{align}
The equation is regular at infinity where $G=0$ if the source term $F$ falls off sufficiently fast towards infinity. The transformed Helmholtz operator is essentially the Helmholtz operator in hyperbolic space \cite{stoll2016harmonic}. 

\subsection{Variable coefficients}

Helmholtz equations with variable coefficients arise in various problems such as seismic full waveform inversion or gravitational waves in curved spacetimes. The Green function may not be readily available in such problems. We can nevertheless apply null infinity compactification if certain asymptotic conditions are fulfilled.  

Consider the one-dimensional case with general coefficients
\[ a(x) \, d_x^2 U + b(x) \, d_x U + c(x) U = 0 \]
For demonstration, we rescale as $u=e^{-\iu h(x)} U$ and map infinity to the origin with $x=1/\rho$. We get 
\[ a \rho^2 d_\rho^2 u + \left( 2 a \rho - b - 2 \iu a H \right) d_\rho u + \left[ \frac{c-a H^2}{\rho^2} + \iu H \frac{b}{\rho^2} + H_\rho a \right] u = 0 \]
The regularity of the equation at infinity requires the coefficients of the Helmholtz equation to be regular at infinity. In addition, we get the fall-off conditions $b=O(x^{-2})$ and $c-a H^2=O(x^{-2})$.

Compactification for equations with variable coefficients is particularly relevant for black hole perturbations where the potential in the equation arising from a non-vanishing spacetime curvature extends to infinity \cite{regge1957stability}. We can nevertheless compute such perturbations at null infinity because the fall-off conditions are satisfied in asymptotically flat spacetimes \cite{ZengFramework, jaramillo2021pseudospectrum}. Similarly, this technique may be helpful when waves propagate in heterogeneous media that fill out the entire solution domain.

\section{A null infinity layer}\label{sec:nil}
\begin{figure}[tbhp]
\centering
\begin{tikzpicture}
\fill [blue!4,even odd rule] (0,0) circle[radius=2cm] circle[radius=1.5cm];
\draw [dashed] (0,0) circle[radius=1.5 cm];
\draw [thick] (0,0) circle[radius=2 cm];
\draw[color=black] (-0.4,0.64) node {$D$};   
\draw[fill=gray!20] plot[smooth cycle,tension=.8]
 coordinates{(-0.8,0.2) (0.6,0.5) (0.5,0) (0.2,-.6)};  
\draw [thick] (0,0) -- (30:2cm) node[anchor=235] {$S$};
\draw [thick] (0,0) -- (60:1.5cm) node[anchor=255] {$R$};
\end{tikzpicture}
\caption{A null infinity layer with a scatterer. The outer boundary at $\rho=S$ corresponds to null infinity. The layer matches to the interior computation at the interface $\rho=R=r$.}
\label{fig:annulus}
\end{figure}
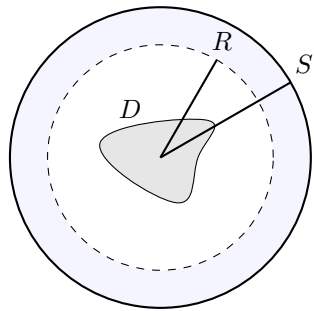

We showed in Sec.~\ref{sec:dispersion} that the transformed wavenumber for incoming waves \eqref{eq:incoming} is higher than the transformed wavenumber for outgoing waves. Near the scatterer, where waves propagate in all directions, null infinity compactification would decrease the accuracy of the numerical solution for incoming waves. Therefore, it is better to use standard coordinates around the scatterer and restrict null infinity compactification to an outer layer (Fig.~\ref{fig:annulus}). Similar thin layers are already commonly used due to the popularity of absorbing and damping layers such as the perfectly matched layer \cite{BERENGER1994185} or the perfect absorbing layer \cite{wang2017perfect, yang2021truly}. In our case, the solution is not artificially damped. Instead, the layer carries outgoing waves to infinity faster than they would otherwise propagate. 

In the time domain, the hyperboloidal layer \cite{ZENGINOGLU20112286} implements this idea. In the frequency domain, we refer to it as a null infinity layer because one can use hyperboloidal and characteristic coordinates by choosing the height function $h$ accordingly. The main feature of the layer is that its outer boundary is at null infinity. 

A null infinity layer solves the Helmholtz equation for a rescaled variable in different coordinates. The solution at the outer boundary of the layer corresponds to the solution at infinity and is therefore of particular interest, especially for radiative problems. We can recover the global solution on the full unbounded domain by reversing the transformations. 

Consider the setup in Fig.~\ref{fig:annulus}. We have a scatterer with boundary $D$. We use standard coordinates in the domain that extends from the scatterer to $r=R$. The boundary of the numerical domain is at $\rho=S$ and the null infinity layer has thickness $S-R$, shaded in light blue. The spatial coordinate mapping in the null infinity layer reads
\begin{align}\label{eq:layer}
g(\rho)=R+ \frac{\rho-R}{\Omega(\rho)} \quad \mathrm{where} \quad \Omega(\rho) = 1 - \frac{(\rho-R)^n}{(S-R)^n}\Theta(\rho-R)\,,
\end{align}
where $\Theta$ is the Heaviside step function and $n$ is an integer power. The height function remains as in \eqref{eq:h} with $K=1$. These choices guarantee a matched interface satisfying
\[ g(R) = R, \quad d_\rho g(R) = 1, \quad d_\rho^n g(R) = 0 \ \mathrm{for} \ n>1. \] 
The optimal choice of the coordinate transformation in the null infinity layer will depend on the problem. For example, \cite{bernuzzi2011binary} uses $n=4$ in the definition of $\Omega$ to get higher derivatives at the layer interface, and \cite{hilditch2018evolution} chooses a smooth transition function instead of the Heaviside function.

\subsection{Relationship to PML and PAL}
The most commonly used absorbing layer in the literature is the perfectly matched layer (PML) \cite{BERENGER1994185}. PML is equivalent to a complex coordinate transformation \cite{chew19943d}. It is very flexible and convenient to use in various coordinates \cite{collino1998perfectly}. Recently, an improved version has been proposed using a complex mapping and rescaling called the perfect absorbing layer (PAL) \cite{wang2017perfect, yang2021truly}. Both PML and PAL damp outgoing waves propagating through the layer.

To contrast these methods with the proposed null infinity layer, we demonstrate the transformations on the example of a one-dimensional, monochromatic, outgoing, plane wave $U(x)=e^{\iu kx}$. We present the calculations for the layer with $x>R$ to avoid carrying Heaviside functions through the expressions. A simple choice of PML reads 
\be\label{eq:pml} x(r) = r + \iu \sigma (r-R), \ee
where $\sigma$ is a free parameter. The PML plane wave solution becomes
\[ U_{\mathrm{PML}}(r) = U(x(r)) = e^{\iu k r} e^{- k \sigma (r-R)}. \]
The wave solution is damped exponentially in the layer where $r>R$. The damping is stronger for thicker layers and larger $\sigma$.

PAL applies a compactification to the complex coordinate transformation of PML (see (3.24) and (3.25) in \cite{yang2021truly}). The mapping is the same as \eqref{eq:layer} with $n=1$ up to a parameter which we set to one 
\[ r(\rho) = R + T(\rho), \quad \mathrm{with} \quad T(\rho) = \frac{(S-R)(\rho-R)}{S-\rho} = \frac{\rho-R}{\Omega(\rho)}\,. \]
The compression mapping leads to singular equations as we discussed in Sec.~\ref{sec:simple}. To remove the infinite oscillations at the domain boundary, PAL introduces a similar rescaling as the null infinity layer
\[ U_{\mathrm{PAL}}(\rho) = e^{-\iu k T(\rho)} U_{\mathrm{PML}}(r(\rho)) = e^{\iu k R} e^{-k \sigma T(\rho)}. \]
The PAL solution is infinitely damped near the outer boundary of the layer irrespective of the thickness, which gives improved decay estimates for the PAL solution than for the PML solution. This property is desirable for constructing a perfect absorbing layer, but it obstructs the recovery of the solution in the exterior domain. We do not need to artificially dampen the outgoing solution when we use a real compactification with the rescaling. The essential benefit of this approach is that the radiative solution becomes directly available at the outer boundary. We can construct the null infinity solution in the layer to take the same form as the solution in standard coordinates using $x=g(\rho)$ and $h(\rho) = g(\rho) - \rho$
\[ u_{\mathrm{NIL}}(\rho) = e^{-\iu k h(\rho)} U(x(\rho)) = e^{-\iu k (g(\rho)-\rho)} e^{\iu k g(\rho)} = e^{\iu k \rho}. \]
We can also recover the PAL solution without the damping by using a slightly different choice for the height function,  $h(\rho) = g(\rho) + R$.  

While the NIL solution looks the same as the standard plane-wave solution, the transformed Helmholtz equation looks different from the Helmholtz equation in standard coordinates. In particular, the transformed equation is regular throughout the domain in compressed coordinates and does not require boundary data. 

\section{Numerical experiments}\label{sec:numerical}
\subsection{Example 1: Plane wave in 1D}\label{sec:oned}
One-dimensional examples are generally not representative for the difficulties related to the outer boundary problem. In our case, however, the essential elements of null infinity compactification are present in one dimension because the transformation acts in the radial direction. Once radial and angular directions are separated, the technique in higher dimensions is similar to its implementation in one dimension because the dynamics of outgoing waves in the asymptotic region is essentially one-dimensional.

Consider the following one-dimensional problem on an unbounded domain
\begin{equation}
	\label{eq:oneD}
	\begin{gathered}
		U_{xx} + k^2 U = 0, \quad \mathrm{in} \ x\in[a,\infty),\\
		U|_{x=a} = \Psi,\\
		\lim_{x\to\infty} (U_x - \iu k U) = 0,
	\end{gathered}
\end{equation}
A simple solution to the above system is the plane wave, $U=e^{\iu k x}$, obtained with the Dirichlet boundary condition $\Psi=e^{\iu k a}$. 

Setting $d=1$ in \eqref{eq:hyp_specific} with vanishing source, we get the transformed problem
\begin{equation}
	\label{eq:hyp_1d}
	\begin{gathered}
		\Omega^2 u_{\rho\rho} - 2 \left(\Omega-\iu k \left(1-\frac{\Omega^2}{K}\right) \right) u_\rho + \left[k^2 \left(\frac{2}{K}-\frac{\Omega^2}{K^2}\right) + 2 \iu k \frac{\Omega}{K} \right] u = 0,\\
		u|_{\rho=\rho_a} = e^{-\iu k (a-\rho_a/K)}\Psi,
	\end{gathered}
\end{equation}
where $\rho_a=\tfrac{a}{1+a}$ and $\Omega(\rho)=1-\rho$. We do not impose the Sommerfeld condition separately because it follows from evaluating the equation at infinity 
\[ u_\rho - \iu \frac{k}{K} u = 0,   \quad \mathrm{at} \ \rho=1. \]
The wavenumber of an outgoing plane wave with null infinity compactification is divided by $K$. The solution reads $u(\rho)=e^{\iu \tfrac{k}{K} \rho}$. We will exploit this modification to compute solutions with high wave numbers.

\begin{figure}[tbhp]
	\centering
	\includegraphics[scale=0.4]{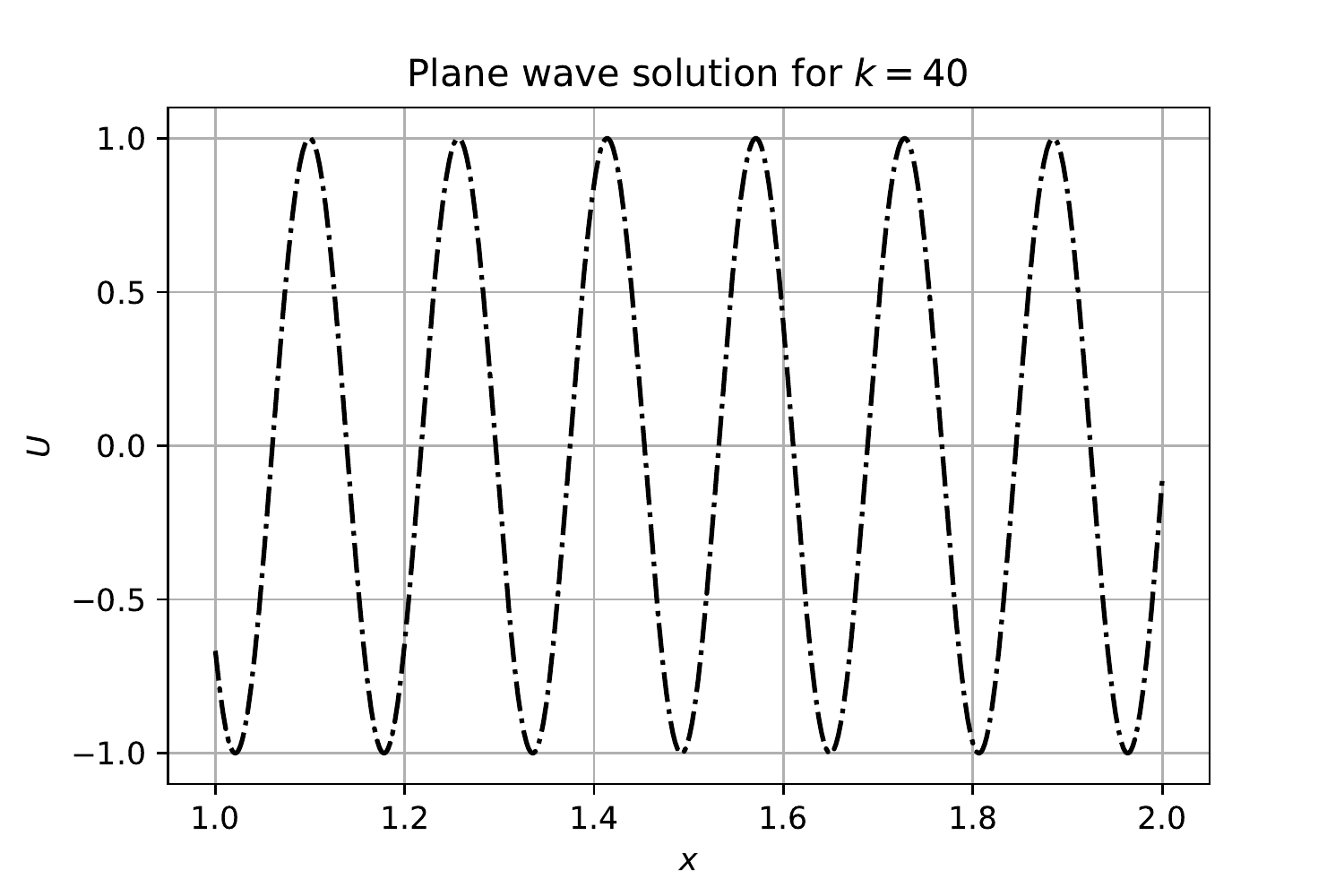}
	\includegraphics[scale=0.4]{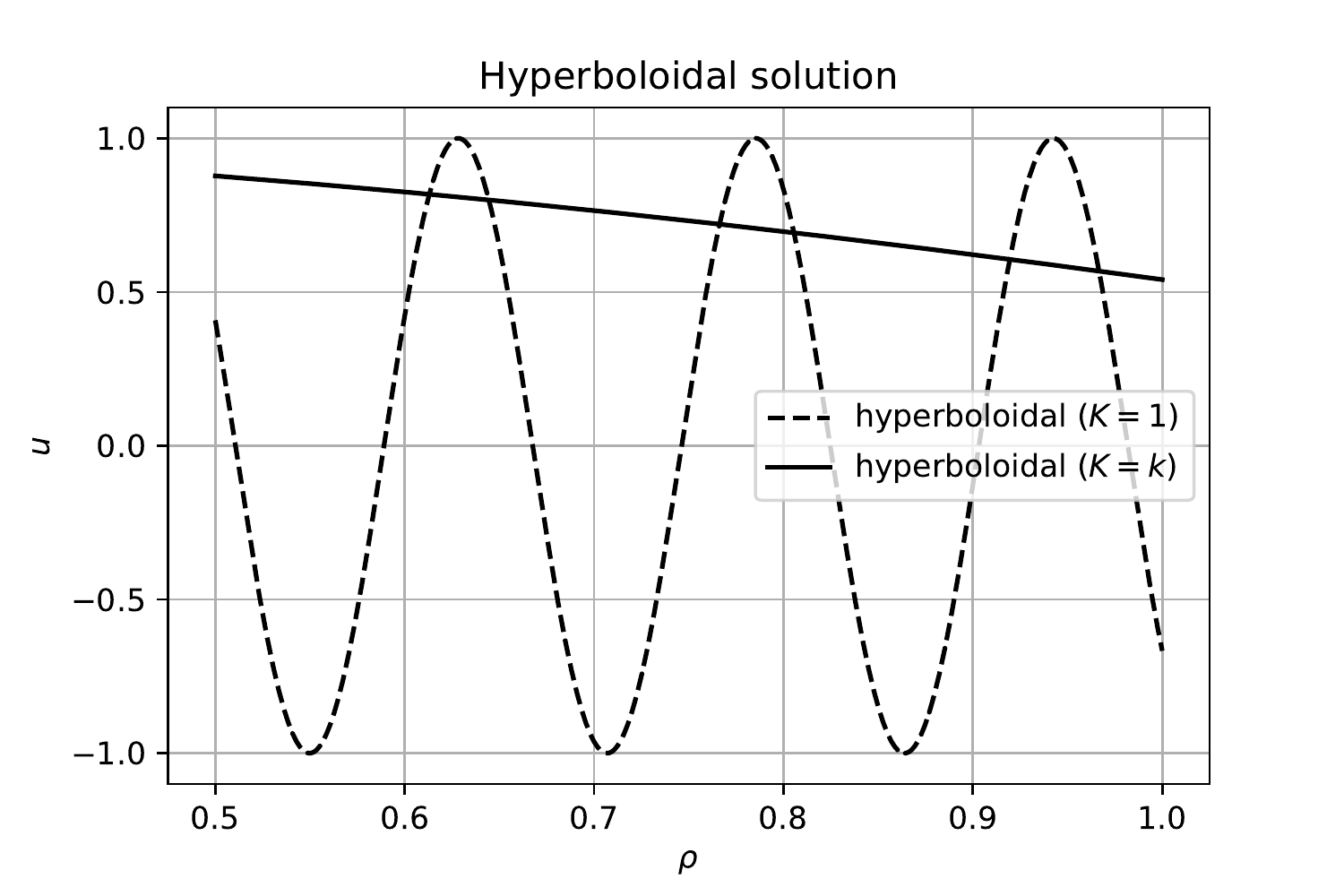}
	\caption{Real parts of the plane wave solution and its transformation in one dimension for $k=40$. The transformed solution is plotted on the right panel for $K=1$ and $K=k$. The domain of the mapped solution corresponds to an infinite domain, and the flattening of the oscillations is controlled by the parameter $K$.}
	\label{fig:oned}
\end{figure}

We plot the real parts of the original solution, $U(x)$, and the transformed solution, $u(\rho)=e^{-ikh(\rho)}U(\rho)$ in Fig.~\ref{fig:oned} with $a=1$ and $k=40$. We restrict the plot in $x$ to $x\in[1,2]$. The domain in $\rho$ includes the point at infinity: $\rho\in[\tfrac{1}{2},1]$.

The plane wave solution has more oscillations across the plotted domain than the hyperboloidal solution, even though the hyperboloidal solution extends to infinity. The freedom in the transformation can be exploited to flatten the oscillations even further as demonstrated by the solid curve on the right panel of Fig.~\ref{fig:oned}, which corresponds to the transformed solution with $K=k$.

We solve \eqref{eq:oneD} and \eqref{eq:hyp_1d} using two numerical discretizations: a finite difference scheme of second order and a spectral-collocation scheme based on Chebyshev polynomials. For solving \eqref{eq:oneD}, the outer boundary data is taken from the exact solution. For solving \eqref{eq:hyp_1d}, the finite difference scheme uses one-sided difference operators at the outer boundary with no boundary data imposed. The spectral-collocation scheme does not require an outer boundary treatment because the principal part is degenerate at infinity which implies a behavioral condition in the terminology of Boyd \cite{boyd2001chebyshev}.

\begin{figure}[tbhp]
	\centering
	\includegraphics[scale=0.4]{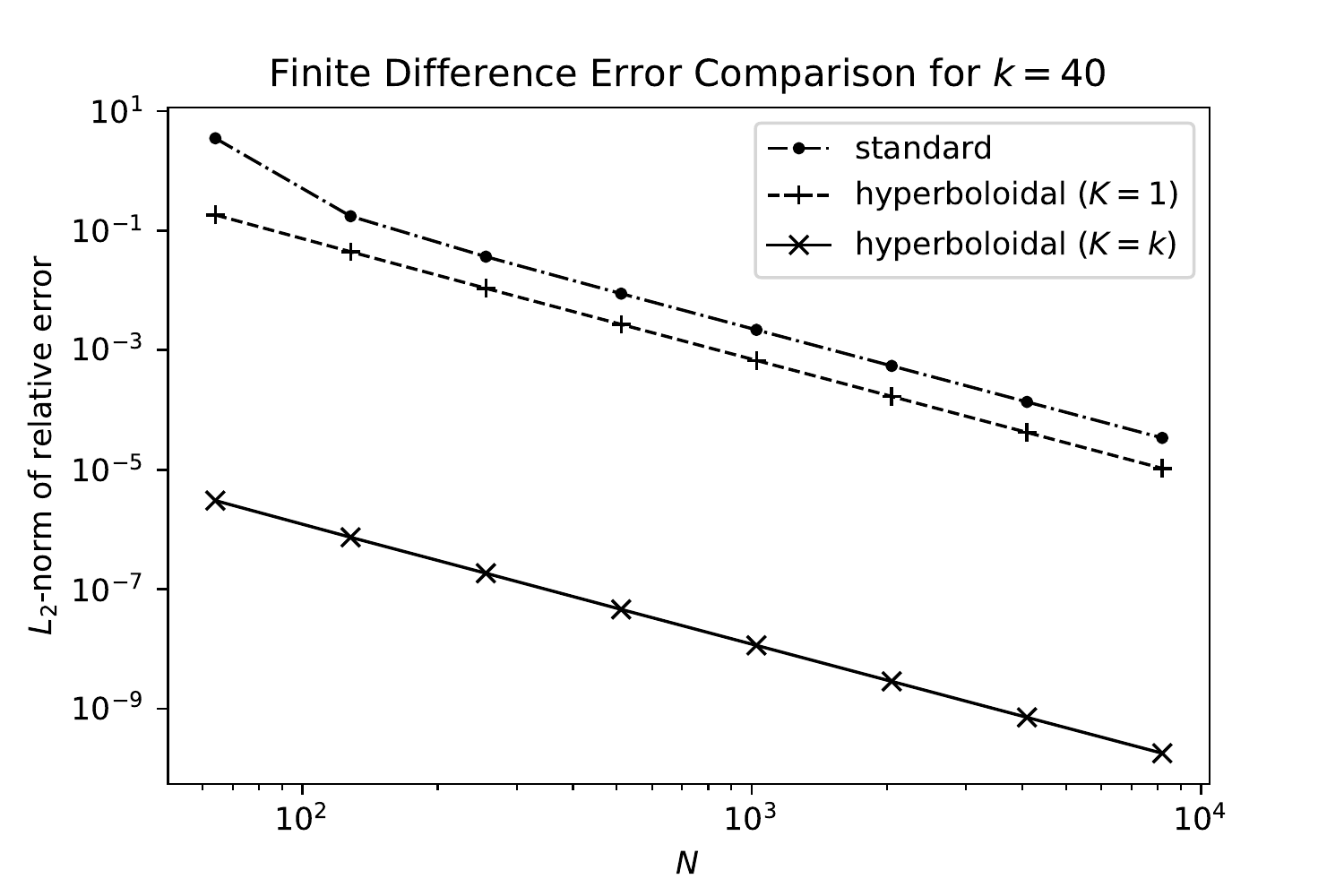}
	\includegraphics[scale=0.4]{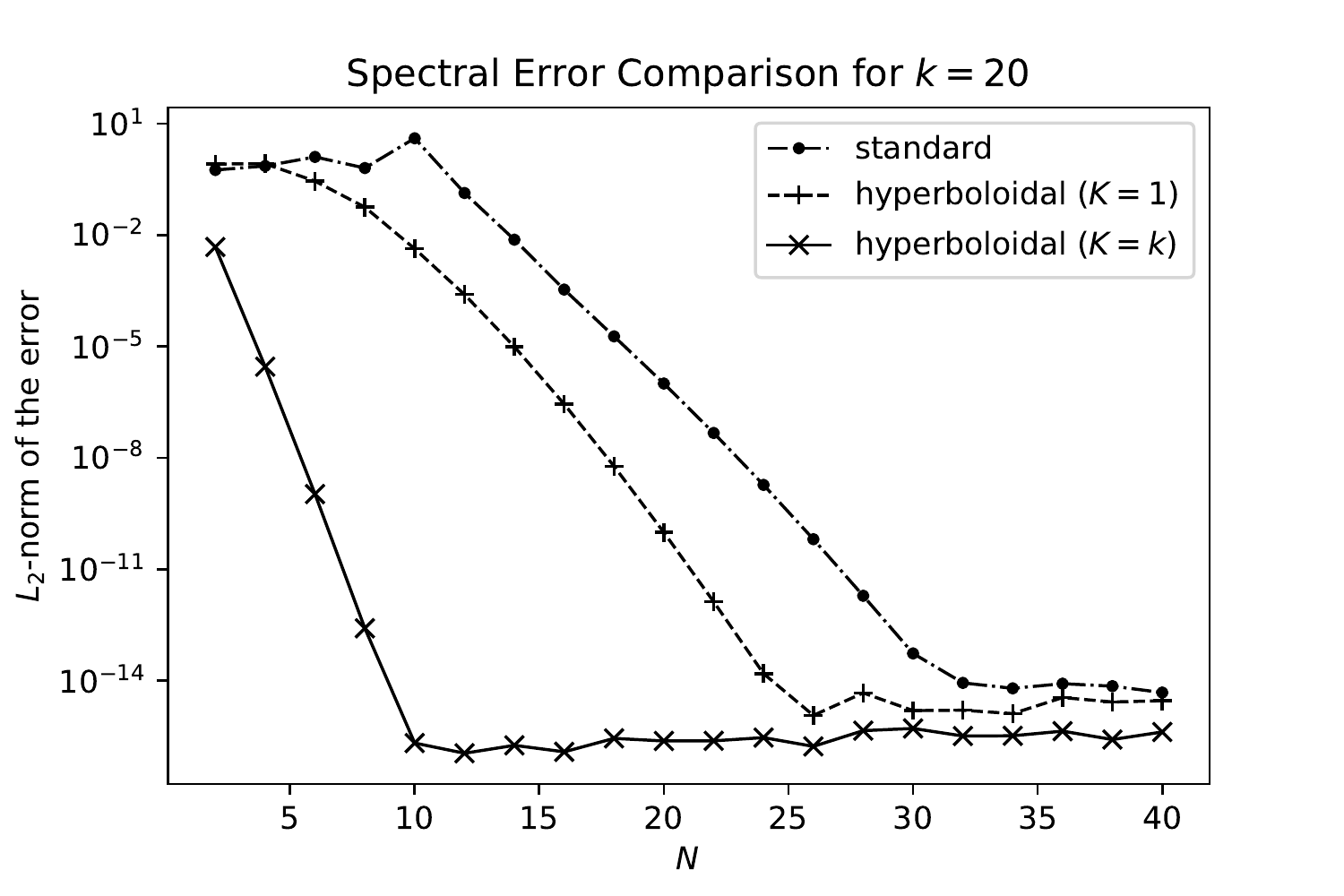}
	\caption{Errors for second order finite difference (left) and Chebyshev spectral (right) solvers in one dimension for $k=40$ for the solution of the 1D problem \eqref{eq:hyp_1d}.}
	\label{fig:errs_oned}
\end{figure}

Figure \ref{fig:errs_oned} shows the corresponding numerical erros for the two schemes. The spectral method is more accurate as expected. The plots also demonstrate the increased efficiency of the transformed solution for $K=k$. In this particular example, one can freely choose a large $K$. In fact, the characteristic foliation with $h(\rho) =g(\rho)$ leads to the transformed solution $u(\rho)=1$. 

\subsubsection{Example 2: Single mode in 2D} 
We consider the source free Helmholtz equation in 2D. We set $d=2$ and $F=0$ in \eqref{eq:helm_sph} and write the equation as a sequence of 1-D equations by expanding the unknown in polar coordinates, $U=\sum U_m(r) e^{\iu m \theta}$. Dropping the subscript $m$, we write
\begin{equation}
	 U_{rr} + U_r + \left(k^2 - \frac{m^2}{r^2}\right)U = 0.
\end{equation}
Applying null infinity compactification using \eqref{eq:g} and \eqref{eq:h} gives 
\begin{equation}
\label{eq:hyp_2d}
	\Omega^2 u_{\rho\rho} - 2 \left(\Omega-\iu k \left(1-\frac{\Omega^2}{K}\right) \right) u_\rho + \left[k^2 \left(\frac{2}{K}-\frac{\Omega^2}{K^2}\right) - \frac{m^2-\tfrac{1}{4}}{\rho^2}  + 2 \iu k \frac{\Omega}{K} \right] u = 0,\\
\end{equation}
Solutions satisfy at the outer boundary the relationship
\[ u_\rho - \iu \left(\frac{k}{K} - \frac{m^2-\tfrac{1}{4}}{2k}\right)u = 0.  \]
This expression suggests that increasing $K$ may not be as effective as in the one dimensional example, especially for large mode numbers. Nevertheless, for high-frequency wave propagation with large $k$, modifying $K$ accordingly should lead to a more efficient solver.

A single mode solution in 2D is given through Hankel functions of the first kind as  $U(r, \theta) = H_m^{(1)}(k r) e^{\iu m \theta}$. After angular decomposition, we get the following expressions for the radial solution in standard and transformed coordinates
\begin{equation}\label{eq:hankel_soln}
	U(r) = H_m^{(1)}(k r), \qquad u(\rho) = \sqrt{r(\rho)} e^{\iu \tfrac{k}{K} \rho} \, \bar{H}_m^{(1)}\left(k r(\rho)\right),
\end{equation}
where $r(\rho)=\rho/\Omega$ and $\bar{H}_{m}^{(1)}$ denotes the exponentially scaled Hankel function defined as $\bar{H}_{m}^{(1)}(z) := e^{-\iu z} H_{m}^{(1)}(z)$. Improved numerical behavior of exponentially scaled Bessel functions for large argument are exploited in many software libraries \cite{amos1986algorithm, 2020SciPy-NMeth}. We can view this exponential scaling geometrically as evaluation along characteristic hypersurfaces of the underlying Minkowski spacetime in line with the time transformation described in step (ii) and the discussion in Sec.~\ref{sec:time} (see also \cite{ZengFramework}).

Many engineering applications require the computation of the far-field pattern. The far-field solution at infinity can be obtained from the asymptotic behavior of the scaled Hankel function. For large $z$, we have \cite{olver1972bessel}
\[ \bar{H}_{m}^{(1)}(z) = \sqrt{\frac{2}{\pi z}} e^{-\iu \pi \left(\tfrac{m}{2}+\tfrac{1}{4} \right)} + O\left(\frac{1}{z^{3/2}}\right) \]
The transformed solution \eqref{eq:hankel_soln} evaluated at null infinity reads
\[ u(1) = \sqrt{\frac{2}{\pi k}} e^{\iu \left(\tfrac{k}{K}- \pi \left(\tfrac{m}{2}+\tfrac{1}{4} \right) \right)}\,. \]
The far-field pattern is directly accessible at the outer boundary of the domain.

\begin{figure}[tbhp]
	\centering
	\includegraphics[scale=0.4]{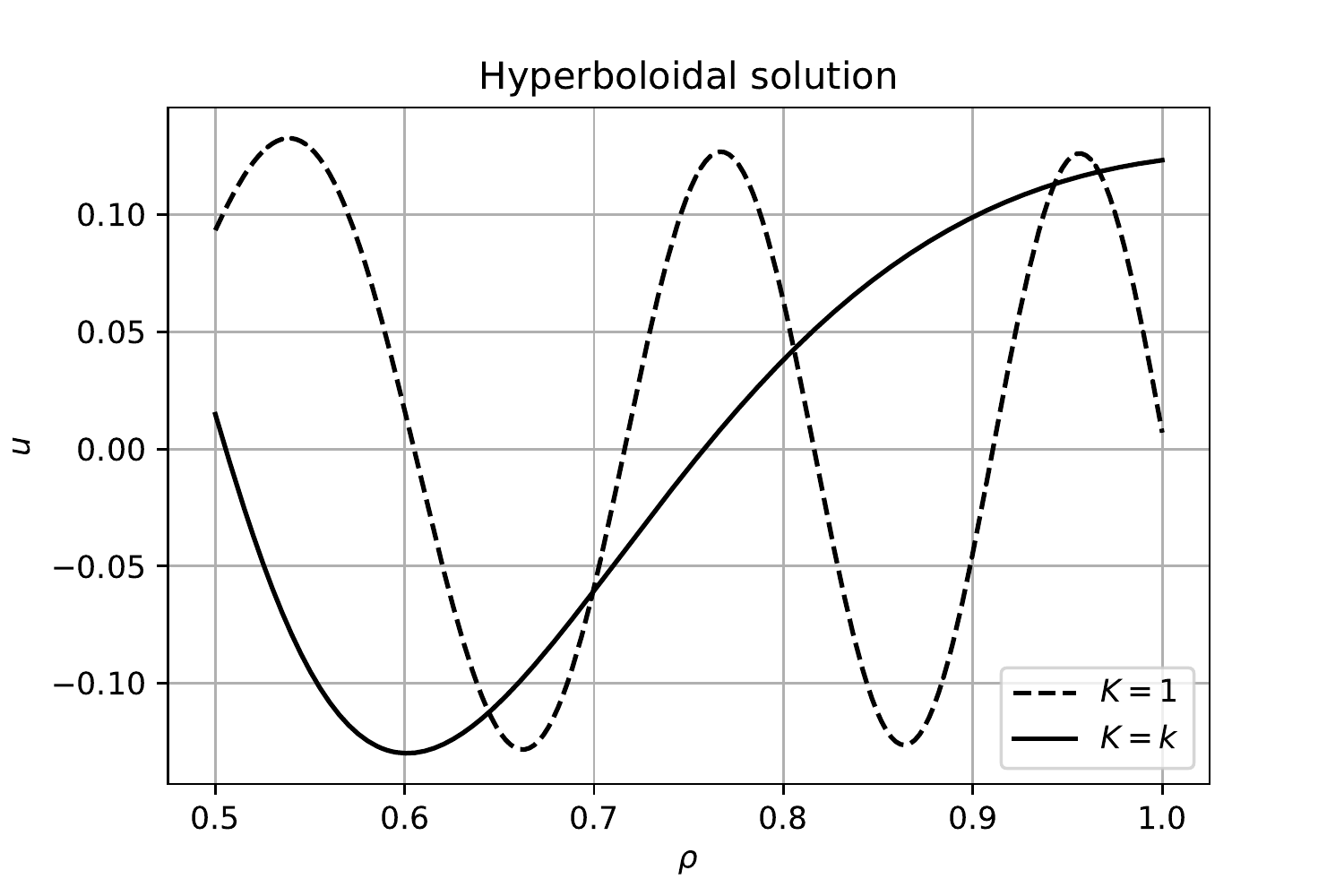}
	\includegraphics[scale=0.4]{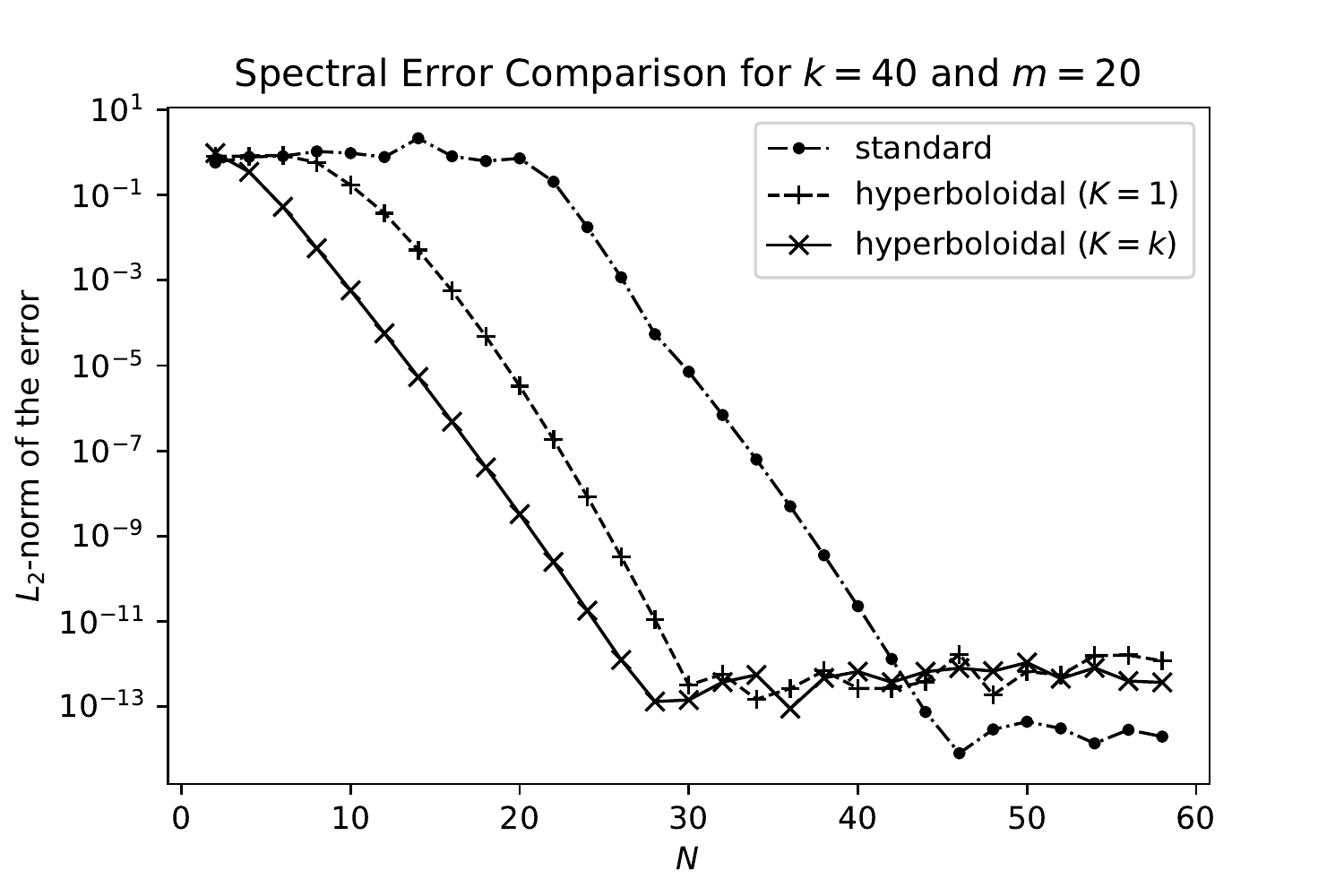}
	\caption{Left panel: Null infinity compactification of the single mode Hankel function for $k=40$ and $m=20$. Due to the mode number $m$, the choice $K=k$ does not remove the oscillations as much as in the 1D case (compare Fig.~\ref{fig:oned}). Right panel: Numerical error comparison between the standard method with exact boundary conditions and the hyperboloidal method with two choices for $K$. The efficiency gain from setting $k=K$ is not as large as in 1D (compare Fig.~\ref{fig:errs_oned}). }
	\label{fig:twod}
\end{figure}

The radial dependence of the transformed solution is plotted for $k=40$ and $m=20$ in Fig.~\ref{fig:twod}. We observe a similar flattening of the oscillations across the domain as in the 1D case. The choice $K=k$ does not remove the oscillations as much as in the 1D case because of the impact of the mode number $m$. The numerical error for the spectral-collocation method is plotted on the right panel of Fig.~\ref{fig:twod}.

\subsubsection{Example 3: Scattering of incoming waves}
In this section, we discuss a scattering problem both in NIC and NIL coordinates. In NIC, the transformation is applied throughout the entire numerical domain while in NIL it is restricted to a thin outer layer.  It may be beneficial in certain applications to use NIL but the layer should be generally avoided if the waves are predominantly outgoing. For the calculations, we set the height function as in \eqref{eq:h} with $K=1$,  and the spatial transformation as in \eqref{eq:g} for NIC and \eqref{eq:layer} for NIL. The equation is given in \eqref{eq:hyp_2d} where we set $K=1$.

\begin{figure}[ht]
  \subfloat[Standard coordinates.]{
	\begin{minipage}[c][1\width]{
	   0.32\textwidth}
	   \centering
	\includegraphics[width=\textwidth]{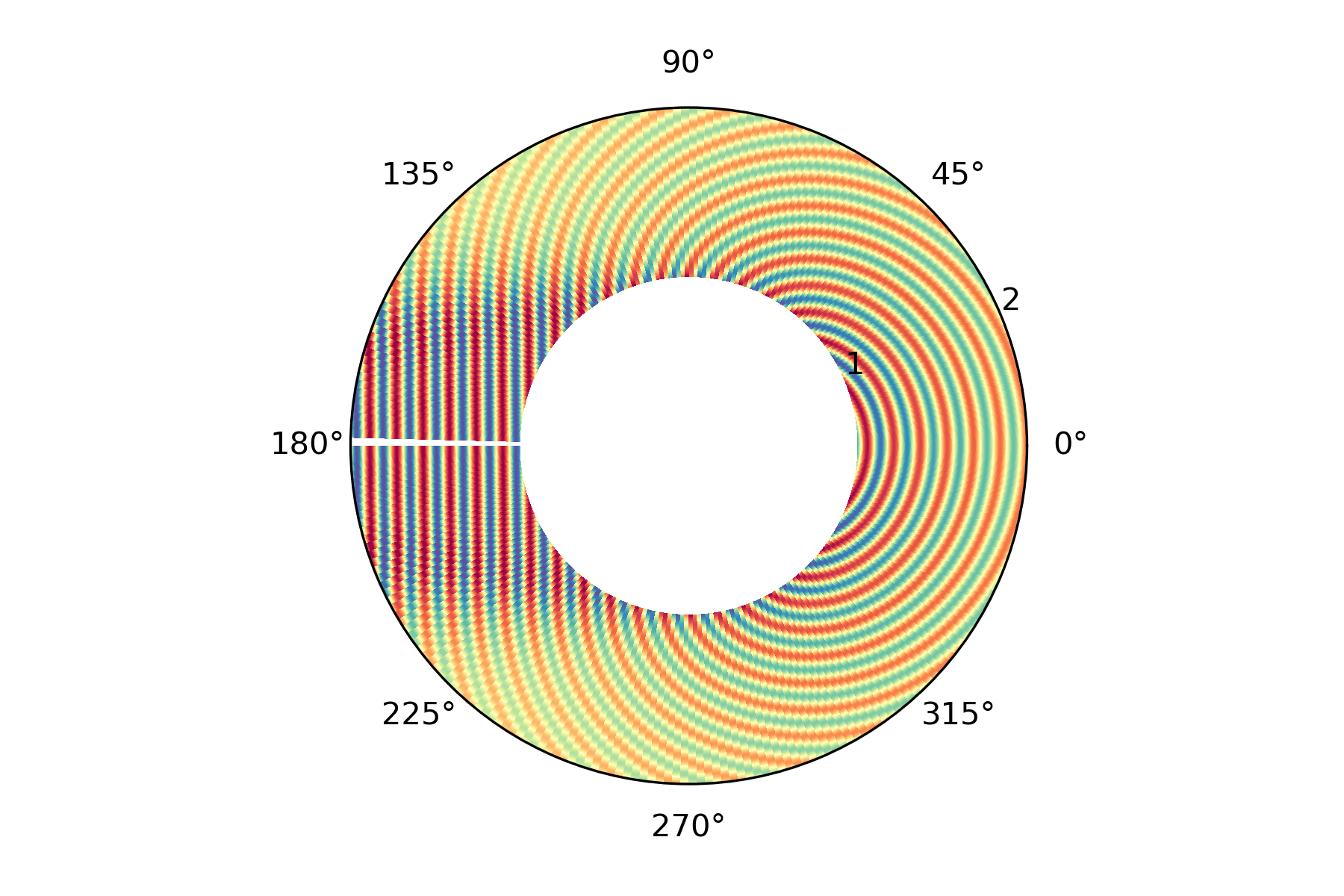}
	\end{minipage}}
 \hfill 	
  \subfloat[NIC coordinates.]{
	\begin{minipage}[c][1\width]{
	   0.32\textwidth}
	   \centering
	   \includegraphics[width=\textwidth]{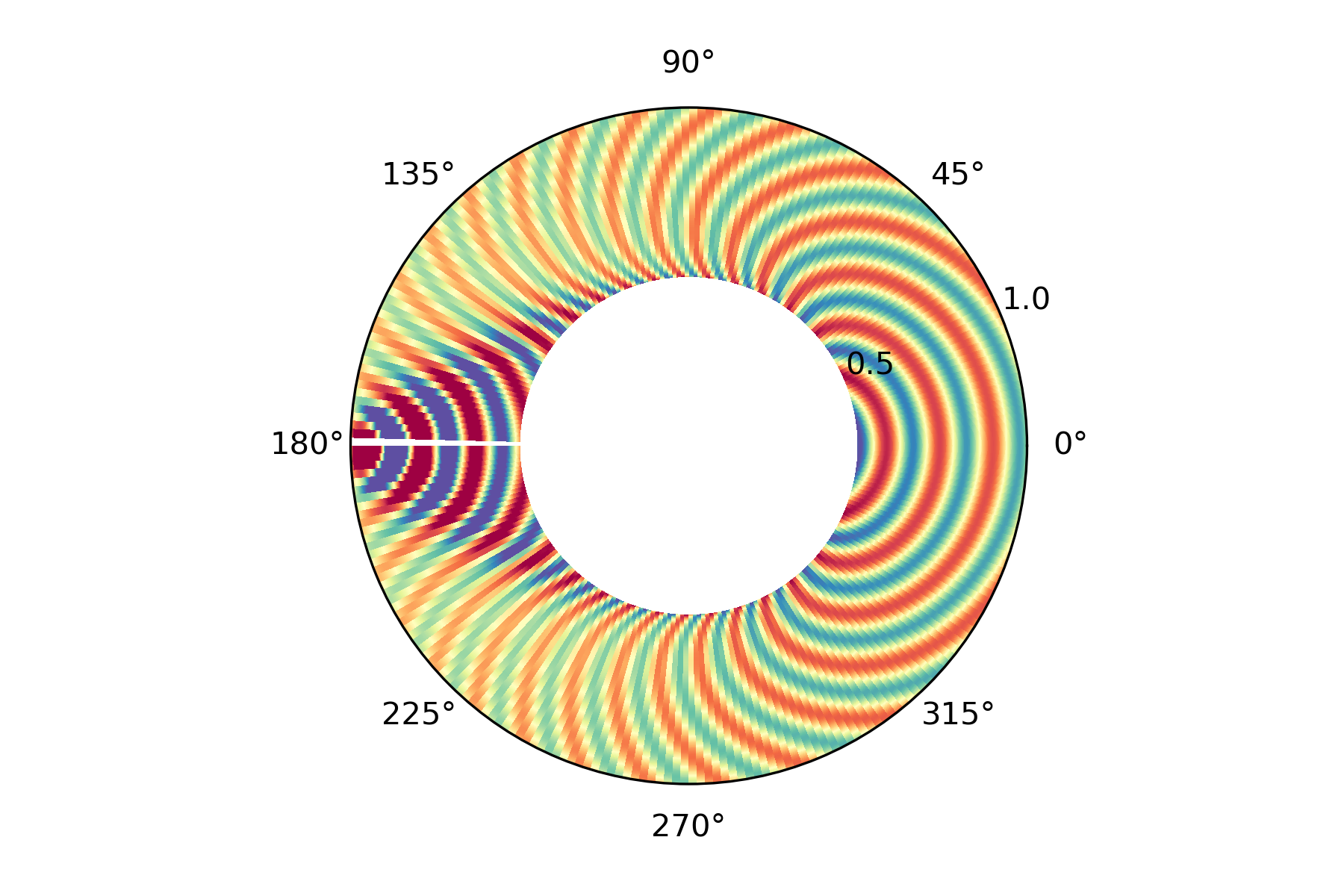}
	\end{minipage}}
 \hfill	
  \subfloat[NIL coordinates.]{
	\begin{minipage}[c][1\width]{
	   0.32\textwidth}
	   \centering
	   \includegraphics[width=\textwidth]{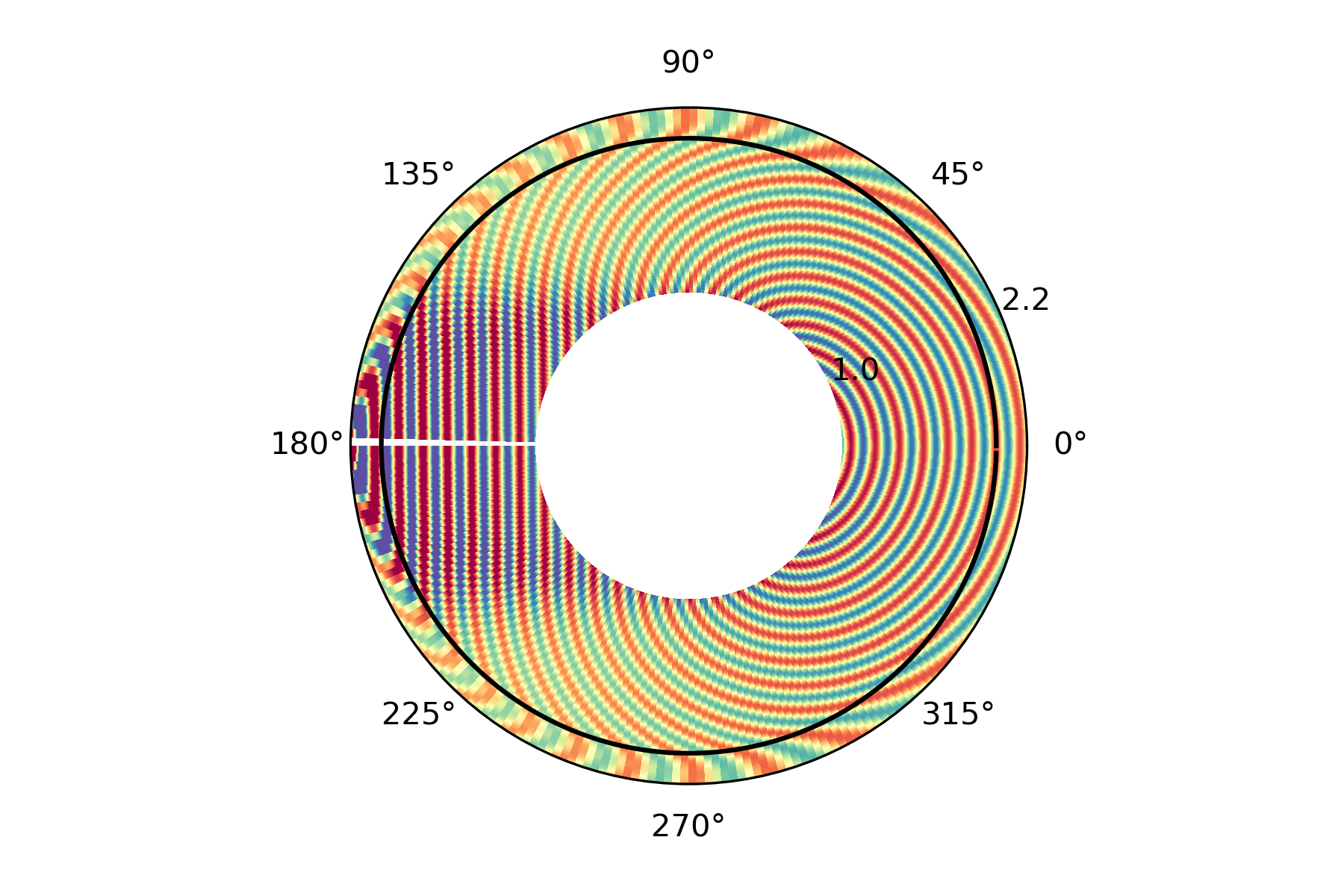}
	\end{minipage}}
\caption{Scattering of an incident wave for $k=40$ in (a) standard coordinates on $r\in[1,2]$; (b) NIC coordinates on $\rho\in[0.5,1]$; and (c) NIL coordinates on $\rho\in[1,2.2]$. The black circle on the NIL solution on panel (c) depicts the interface at $R=2$. The solution for $\rho<R$ matches identically the solution depicted on panel (a). The harmonic content of the NIC representation on (b) is smaller than the NIL representation on (c), which leads to a more efficient computation.}
\label{fig:soft_sound}
\end{figure}

We consider a plane wave in the $x$-direction, $e^{\iu k x} = e^{\iu k r \cos\theta}$, scattered off a circle of radius $R_0=1$. The incident plane wave does not satisfy the Sommerfeld radiation condition in 2D and therefore is not regular at null infinity. The scattered solution is outgoing at infinity, and is therefore regular. The total field is the sum of the impinging field and the scattered field, implying the following Dirichlet boundary condition for the scattered field:
\begin{equation}\label{eq:dirichlet}
	U(r,\theta)|_{r=R_0} = - e^{\iu k R_0 \cos\theta}.
\end{equation}
We write the solution satisfying the Sommerfeld radiation condition as a series expansion (as in, for example, \cite{britt2010compact, yang2021truly})
\[ U(r,\theta) =  \sum_{|m|=0}^{\infty} c_m H_{m}^{(1)}(k r) e^{\iu m \theta} \,, \]
where the coefficients, $c_m$, are determined from a Fourier expansion by requiring that the Dirichlet condition \eqref{eq:dirichlet} is satisfied at the scattering surface:
\[ c_m = \frac{1}{H_{m}^{(1)}(k R_0)} \frac{1}{2\pi} \int_{-\pi}^{\pi} -e^{\iu k R_0 \cos\theta} e^{-\iu m \theta} d\theta = 
- \frac{\iu^m J_m(k R_0)}{H_{j}^{(1)}(k R_0)}. \]
The $J_m$ are Bessel functions. We obtain the transformed initial data and solution by scaling-out the oscillatory decay and compactifying as in \eqref{eq:hankel_soln}. The initial data reads
\[ u(\rho_0,\theta)= \sqrt{R_0} e^{-\iu k (R_0-\rho_0)} U(R_0,\theta) 
\]
The transformed solution reads
\begin{equation}\label{eq:soft_soln}
	u(\rho,\theta) = \sqrt{r(\rho)} e^{\iu k\rho} \sum_{|m|=0}^{\infty} c_m \bar{H}_{m}^{(1)}\left(k r(\rho)\right) e^{\iu m \theta},
\end{equation}
Three representations of the solution in (a) standard, (b) NIC, and (c) NIL coordinates are plotted in Fig.~\ref{fig:soft_sound}. The black interior circle on panel (c) marks the interface at $\rho=R$. The solution for $\rho<R$ agrees identically with the standard representation on panel (a). 

\vspace{-1cm}

\begin{figure}[tbhp]
	\centering
	\includegraphics[scale=0.41]{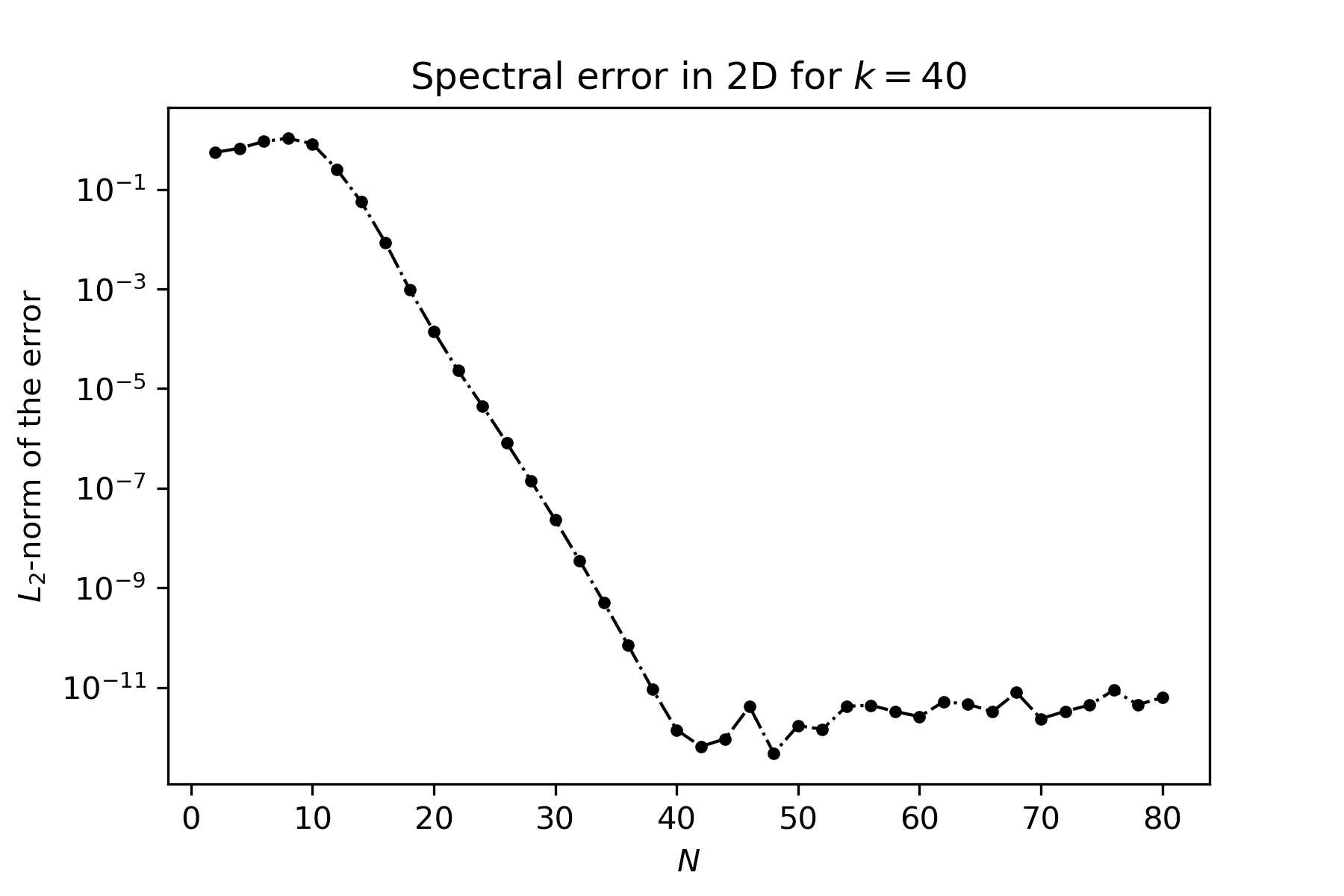}
	\includegraphics[scale=0.41]{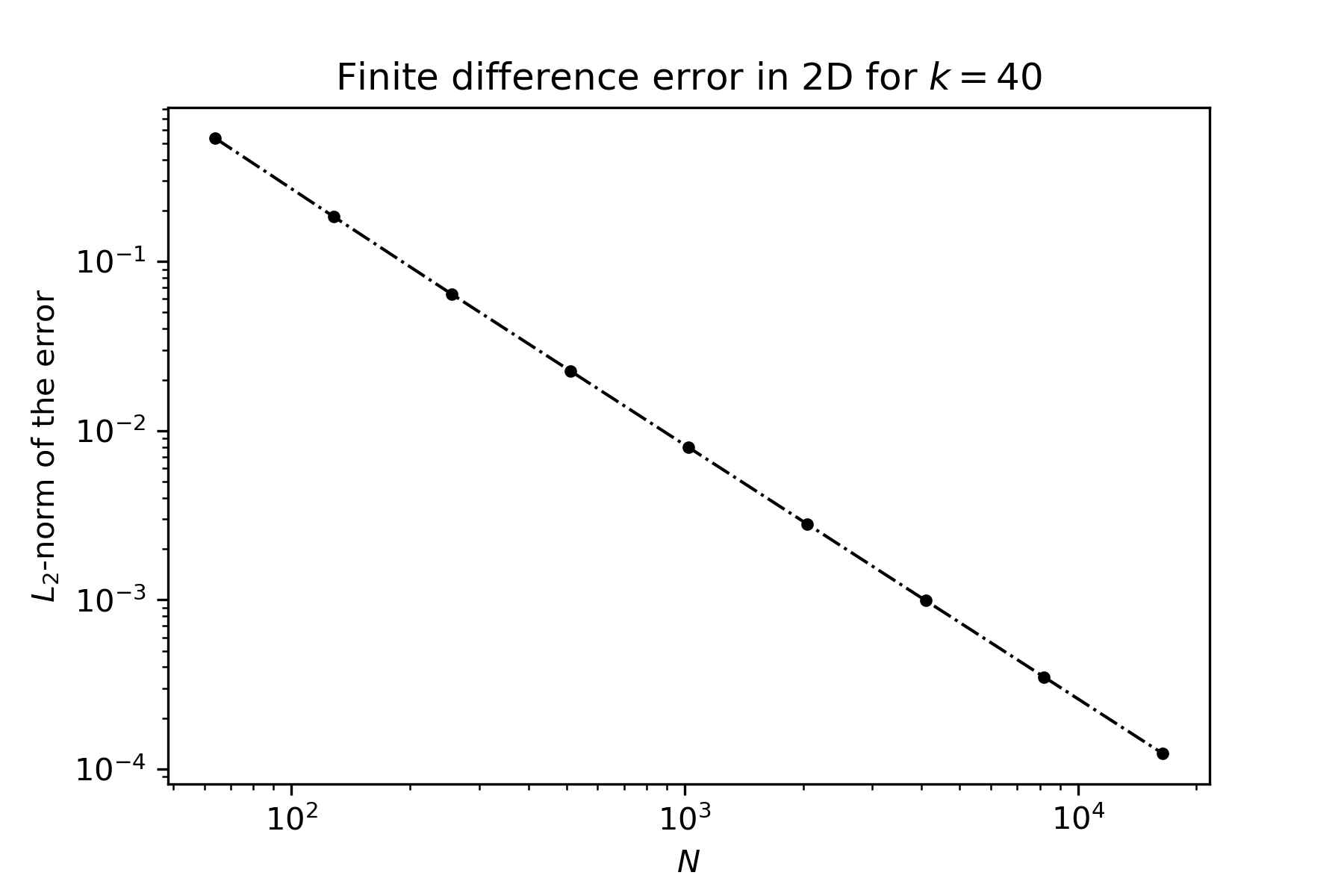}
	\caption{Convergence for the solution of the scattering of an incident plane wave with $k=40$ using spectral (left) and 2nd order finite difference (right) methods. }
	\label{fig:twod_scattering}
\end{figure}

The NIC representation shows fewer waves along the radial directions in accordance with the flattening property of hyperboloidal compactification. The NIL representation distorts the waves in the layer transporting them to null infinity. It also has a higher harmonic content than the NIC representation. Therefore, it is preferable for this particular problem to use the NIC representation. We present in Fig.~\ref{fig:twod_scattering} convergence results for solving the Helmholtz equation in NIC representation both with spectral (left panel) and 2nd order finite difference (right panel) methods. Spectral methods are clearly superior in such problems.

\section{Summary and Discussion}
The main idea of this paper is the following directive: Scale-out the oscillatory decay and compactify \eqref{eq:summary}. Applying this simple directive leads to equation \eqref{eq:hyp_general} which has various advantages for scientific computation. First of all, the equation does not require the formulation of an artificial outer boundary problem. One solves the degenerate equation without imposing boundary data because the equation geometrically incorporates the no-incoming radiation condition through a behavioral boundary. This property of the equation also simplifies its discretization. For example, among the motivations for constructing compact finite difference operators with narrow stencils is the boundary treatment \cite{britt2010compact}. In contrast, null infinity compactification allows us to apply high-order finite-difference operators using one-sided stencils near the outer boundary. In spectral methods, the function space automatically incorporates the behavior near the boundary.

Another advantage of null infinity compactification is access to the asymptotic solution. In many scientific and engineering applications, we need to measure the outgoing radiation, such as the echo area in acoustics, the radar cross-section in electromagnetics, or the gravitational waveform in general relativity. Truncating the problem domain makes the extraction of such quantities cumbersome. With null infinity compactification, we can read off the radiation from the numerical solution performed on a relatively small domain directly without post-processing. 

The flexibility of the method allows us to adapt the coordinates to the problem. For example, black hole perturbations propagate both to infinity and the black hole. We can adjust the height function to incorporate this behavior into the equation resulting in behavioral boundaries at both ends of the domain (black hole horizon and null infinity)\cite{ZengFramework}. Similarly, suppose the geometry of a wave scattering problem has preferred directions of wave propagation. In that case, one can incorporate these into the height function to improve the efficiency of the numerical method.


The relationship of rescalings with time transformations puts existing literature on Helmholtz equations into a geometric framework. We understand the exponential rescaling in PAL and infinite elements as transformations to a characteristic coordinate system. The rescaling of the asymptotic decay corresponds to the conformal compactification of the underlying spacetime. This unified framework allows us to extend existing methods to general height functions to incorporate variable propagation speeds and heterogeneous media.
 
Currently, the method is limited to smooth outer boundaries. This restriction comes from the topology of null infinity: cross-sections of future null infinity have spherical topology \cite{Penrose65, Geroch}. To simulate scattering from long objects, we can use ellipsoidal or prolate spheroidal coordinates. More generally, we can parametrize a smooth boundary as demonstrated in \cite{yang2021truly}. However, null infinity compactification as presented in this paper is not readily applicable to Cartesian or polygonal domains.

The method should be further investigated in real-world applications where restriction of the transformations to a layer may be necessary. Another avenue of research is the numerical analysis of the transformed equation \eqref{eq:hyp_general} for suitable choices of height function and compactification. The nature of the transformations and the outcomes of numerical experiments suggest that null infinity compactification is local, exact, and stable. A detailed numerical analysis of these properties is lacking. Such research would also be helpful to obtain guidelines for choices of free functions and parameters involved in the transformations, such as the thickness of the null infinity layer or the shape of the height function.

\vspace{-7mm}

\section*{Acknowledgments}
I thank the Institute for Physical Science and Technology and its director Konstantina Trivisa for providing me space and time to perform this research. I thank Rodrigo Panosso Macedo, Ricardo Nochetto, and Eitan Tadmor for useful discussions, and Jie Shen for reminding me of the references \cite{wang2017perfect, yang2021truly}. 

\pagebreak
\bibliographystyle{siam}

\bibliography{refs}

\end{document}